\documentclass{amsproc}
\usepackage{amsmath,amssymb}

\newtheorem{theorem}{Theorem}[section]

\newtheorem{proposition}[theorem]{Proposition}
\newtheorem{corollary}[theorem]{Corollary}

\theoremstyle{remark}
\newtheorem{remark}[theorem]{Remark}

\numberwithin{equation}{section}

\newcommand {\bbR} {\mathbb{R}}

\newcommand {\calD} {\mathcal{D}}
\newcommand {\calM} {\mathcal{M}}
\newcommand {\calS} {\mathcal{S}}

\renewcommand {\Re}{\mathop{\mathrm{Re}}}

\begin{document}

\title[Zeroes of $L$-functions of modular forms]{Uniform distribution of zeroes of $L$-functions of modular forms}
\author{Alexey Zykin}
\address{
Alexey Zykin
\newline \indent
Laboratoire GAATI, Universit\'e de la Polyn\'esie fran\c caise, BP 6570 --- 98702 Faa'a, Tahiti, Polyn\'esie fran\c caise
\newline \indent
Department of Mathematics of the National Research University Higher School of Economics
\newline \indent
AG Laboratory NRU HSE 
\newline \indent
Laboratoire Poncelet (UMI 2615)
\newline \indent
Institute for Information Transmission Problems of the Russian Academy of Sciences
}

\email{alzykin@gmail.com}
\thanks{The author was partially supported by AG Laboratory NRU HSE, RF government grant, ag.  11.G34.31.0023, by ANR Globes ANR-12-JS01-0007-01.}
\date{}

\begin{abstract}
We prove under GRH that zeros of $L$-functions of modular forms of level $N$ and weight $k$ become uniformly distributed on the critical line when $N+k\to\infty.$
\end{abstract}

\maketitle
\section{Introduction}
It is well known that zeroes of $L$-functions contain an important information about the arithmetic properties of the objects to which these $L$-functions are associated. The question about the distribution of these zeroes on the critical line was studied by many authors. This problem can be looked upon from many angles (the proportion of zeroes on the critical line, low zeroes, zero spacing, etc.).

In this paper we study the distribution of zeroes of $L$-functions on the critical line when we let vary the modular form to which the $L$-function is associated. The same question was considered by S. Lang in \cite{Lan71} and M. Tsfasman and S. Vl\u adu\c t in \cite{TV} for the Dedekind zeta function of number fields.

Let $f(z)$ be a holomorphic cusp of weight $k=k_f$ for the group $\Gamma_0(N)$ such that $f(z)=\sum\limits_{n=1}^{\infty}a_n n^{(k-1)/2} e^{2\pi i n z}$  is its normalized Fourier expansion at the cusp $\infty.$ We suppose that $f(z)$ is a primitive form in the sense of Atkin--Lehner \cite{AL} (it is a new form and a normalized eigen form for all Hecke operators), so $L_f(s)$ can be defined by the Euler product
$$L_f(s)=\prod\limits_{p\mid N}(1-a_p p^{-s})^{-1} \prod\limits_{p\,\nmid N}(1-a_p p^{-s}+p^{-2s})^{-1}.$$
We denote by $\alpha_p$ and $\bar{\alpha}_p$ the two conjugate roots of the polynomial $1-a_p p^{-s}+p^{-2s}$. Deligne has shown (see \cite{Del}) that $|\alpha_p|=|\bar{\alpha}_p|=1$ for $p \nmid N$ (the Ramanujan--Peterson conjecture). On the other hand, one knows (see \cite{AL}) that for $p \mid N$ we have $|a_p|\leq 1.$

If we define the gamma factor by
$$\gamma_f(s)=\pi^{-s}\Gamma\left(\frac{s+(k-1)/2}{2}\right)\Gamma\left(\frac{s+(k+1)/2}{2}\right)=c_k(2\pi)^{-s}\Gamma\left(s+\frac{k-1}{2}\right)$$
with $c_k=2^{(3-k)/2}\sqrt{\pi},$ then the function $\Lambda(s)=N^{s/2}\gamma_f(s)L_f(s)$ is entire and satisfies the functional equation $\Lambda(s)=w\Lambda(1-s)$ with $w=\pm 1.$ The Generalized Riemann Hypothesis (GRH) for $L$-function of modular forms states that all the non-trivial zeroes of these $L$-functions lie on the critical line $\Re s = \frac{1}{2}.$ Throughout the paper we assume that GRH is true.

The analytic conductor $q_f$ (see \cite{IK}) is defined as
$$q_f=N\left(\frac{k-1}{2}+3\right)\left(\frac{k+1}{2}+3\right)\sim \frac{Nk^2}{4},$$
when $k\to \infty.$ We will use the last expression (or, more precisely, its logarithm minus a constant) as a weight in all the zero sums in the paper.

To each $f(z)$ we can associate the measure
$$\Delta_f:=\frac{2\pi}{\log q_f}\sum_{L_f(\rho)=0} \delta_{t(\rho)},$$
where $t(\rho)=\frac{1}{i}\left(\rho-\frac{1}{2}\right)$ and $\rho$ runs through all non-trivial zeroes of $L_f(s)$; here $\delta_a$ denotes the atomic (Dirac) measure at $a.$ Since we suppose that GRH is true, $\Delta_f$ is a discrete measure on $\bbR.$ Moreover, it can easily be seen that $\Delta_f$ is a measure of slow growth (see below).

Our main result is the following one:

\begin{theorem}
\label{ZEROESMODULAR}

Assuming GRH, for any family $\{f_j(z)\}$ of primitive forms with $q_{f_j}\to \infty$ the limit
$$\Delta = \lim_{j\to \infty} \Delta_j= \lim_{j\to \infty} \Delta_{f_j}$$
exists in the space of measures of slow growth on $\bbR$ and is equal to the measure with density $1$ (i. e. $dx$).
\end{theorem}

\section{Proof of theorem \ref{ZEROESMODULAR}}
Our method of the proof will, roughly speaking, follow that of \cite{TV}, where a similar question is treated in the case of Dedekind zeta functions. It will even be simplier in our case due to the fact that the family we consider is "asymptotically bad".

Let us recall a few facts and definitions from the theory of distribution. We will use \cite{Sch} as our main reference. Recall that the Schwartz space $\calS=\calS(\bbR)$ is the space of all real valued infinitely differentiable rapidly decreasing functions on $\bbR$ (i.~e. $\phi(x)$ and any its derivative go to 0 when $|x|\to \infty$ faster then any power of $|x|$). The space $\calD(\bbR)$ is defined to be the space of all real valued infinitely differentiable functions with compact support on $\bbR.$ Both $\calS(\bbR)$ and $\calD(\bbR)$ are equipped with the structures of topological vector spaces.

The space $\calD'$ (resp. $\calS'$), topologically dual to $\calD$ (resp. $\calS$) is called the space of distribution (resp. tempered distributions). We also define the space of measures $\calM$ as the topological dual of the space of real valued continuous functions with compact support on $\bbR.$ The space $\calM$ contains a cone of positive measures $\calM_{+},$ i.~e. of measures taking positive values on positive functions. One has the following inclusions: $\calS'\subset\calD'$ and $\calM_{+} \subset \calM \subset\calD'.$ The intersection $\calM_{sl}=\calM\cap\calS'$ is called the space of measures of slow growth. A measure $\mu$ of slow growth can be characterized by the property that for some positive integer $k$ the integral
$$\int_{-\infty}^{+\infty}(x^2+1)^{-k} d\mu$$
converges (see \cite[Thm. VII of Ch. VII]{Sch}). In particular, from this criterion and the fact that the series $\sum\limits_{\rho\neq 0, 1}|\rho|^{-2}$ converges (\cite[Lemma 5.5]{IK}), we see that $\Delta_f$ is a measure of slow growth for any $f.$

Finally, we note that the Fourier transform $\hat{}$ is defined on $\calS$ and $\calS'$ and is a topological automorphism on these spaces. $\calD$ is known to be dense in $\calS$ and so $\hat{\calD}$ is also dense in $\calS=\hat{\calS}.$ To check that $\mu$ is a measure of slow growth it is enough to check that it is defined on a dense subset and that it is continuous on this dense subset in the topology of $\calS.$ In the same way, to check that a sequence of measures of slow growth converges to a measure of slow growth it is enough to check its convergence on a dense subset to a measure continuous on this dense subset. This follows from the definition of measures as linear functionals.

Our main tool will be a version of Weil explicit formula for $L$-functions of modular forms proven in \cite[I.2]{Mes} or in \cite[theorem 5.12]{IK} (in the last source some extra conditions on test functions are imposed).

Suppose $F\in \calS(\bbR)$ satisfies for some $\epsilon > 0$ the following condition
\begin{equation}
\label{cond}
|F(x)|, |F'(x)| \ll c e^{(-\frac{1}{2}+\epsilon)|x|} \text{ as } |x|\to\infty.
\end{equation}
Let
$$\Phi(s):=\int_0^{\infty} F(x)e^{(s-\frac{1}{2})x} dx = \hat{F}(t),$$
where $s=\frac{1}{2}+it.$ The next proposition gives us the explicit formula that we need to relate the sum over zeroes to the sum of coefficient of modular forms:

\begin{proposition}
\label{propexpl}
Let $f(z)$ be a primitive form of level $N$ and weight $k.$ Then the limit
$$\sum_{L_f(\rho)=0}\Phi(\rho)=\lim_{T\to \infty}\sum_{\substack{L_f(\rho)=0\\ |\rho|< T}}\Phi(\rho)$$
exists and we have the following formula:
\begin{multline*}
\sum_{L_f(\rho)=0}\Phi(\rho) = -\sum_{p,m}b(p^m)(F(m\log p)+F(-m\log p))\frac{\log p}{p^{m/2}}+\\
+F(0)(\log N - 2 \log (2\pi))+\frac{1}{\pi} \int_{-\infty}^{+\infty} \frac{\Phi\left(\frac{1}{2}+it\right)+\Phi\left(\frac{1}{2}-it\right)}{2}\cdot\psi\left(\frac{k}{2}+it\right)dt,
\end{multline*}
where $\psi(s)=\Gamma'(s)/\Gamma(s)$, $b(p^m)=(a_p)^m$ if $p \mid N$ and $b(p^m)=(\alpha_p)^m+(\bar{\alpha}_p)^m$ otherwise.
\end{proposition}

Taking a subsequence of $\{f_j\}$ we can assume that the limit $\alpha=\lim\limits_{j\to\infty}\frac{\log N_j}{\log N_j+\log k_j}$ exists. We will check the convergence of measures on $\hat{\calD}.$ From the above discussion this is enough to prove the result. Let us take any $\phi\in \hat{\calD},$  $\phi=\hat{F},$ $F \in \calD.$ We have $\phi(t)=\Phi\left(\frac{1}{2}+it\right).$ The function $F$ satisfies the condition (\ref{cond}), so we can apply the explicit formula to it. We fix $\phi(t)$ and let vary $f_j$ Then, we get the equality when $j\to \infty$.
\begin{equation}
\label{eqlim}
\Delta(\phi)= 2\pi F(0) \alpha + 2\int_{-\infty}^{+\infty} \frac{\phi(t)+\phi(-t)}{2}\cdot\lim_{j\to\infty}\frac{\psi\left(\frac{k_j}{2}+it\right)}{\log N_j+\log k_j}dt,
\end{equation}
since $|b(p^m)|\leq 2$ and the integral is uniformly convergent as $\phi(t)\in \calS.$ The limit under the integral sign can be evaluated using the Stirling formula $\psi(s)=\log s+O\left(\frac{1}{|s|}\right)$ (see \cite[p. 332]{Lan94}). This gives us
$$\lim_{j\to\infty}\frac{\psi\left(\frac{k_j}{2}+it\right)}{\log N_j+\log k_j}=\frac{1}{2}(1-\alpha).$$
But $\int_{-\infty}^{+\infty}\psi(t) dt = 2\pi F(0)$ and so the right hand side of (\ref{eqlim}) equals
$$2\pi F(0)\alpha+ 2\pi F(0)(1-\alpha)=2\pi F(0)=\int_{-\infty}^{+\infty}\phi(t) dt.$$
This concludes the proof of the theorem. \qed

\begin{corollary}
Any fixed interval around $s=\frac{1}{2}$ contains zeroes of $L_{f}(s)$ if $q_f$ is sufficiently large.
\end{corollary}

\begin{remark}
One can prove a similar equidistribution statement for $L$-functions of bounded degree in the Selberg class, assuming suitable conjectures (like the Generalized Riemann Hypothesis). It is an interesting question how zeroes of $L$-functions are distributed if the degree of these $L$-functions grows with the analytic conductor. Some examples of non-trivial distributions of zeroes for Dedekind zeta functions are considered in \cite{TV}.
\end{remark}
%We are left to show that For any $a\in \bbR$ and $y \geq 0$ let us define the function $H_{y,a}(x)\in\calS(\bbR)$  by
%$$H_{y,a}(x):=\frac{1}{2\sqrt{\pi y}}\exp\left(\frac{-(a-x)^2}{4y}\right).$$

%We have the following lemma (see \cite[Lemma 5.2]{TV}):
%\begin{lemma}
%Let $\mu$ be a positive measure on $\bbR$ such that for any $a\in\bbR$ one has
%$$\lim_{y\to+0}\mu(H_{y,a})=M(a)$$
%for some continuous function $M:\bbR\to\bbR.$ Then $\mu=M(x)dx.$
%\end{lemma}
\bibliographystyle{amsalpha}

\end{document}